\documentclass[12pt,draftcls,onecolumn]{article}

\usepackage{mathrsfs}
\usepackage{amsfonts}
\usepackage{amssymb}
\usepackage{amssymb,amsmath,amsthm, amsfonts}
\usepackage{color}
\usepackage{graphicx}

\def\sqr#1#2{{\vcenter{\vbox{\hrule height.#2pt
              \hbox{\vrule width.#2pt height#1pt \kern#1pt \vrule width.#2pt}
          \hrule height.#2pt}}}}

\usepackage{amsmath,amssymb,amsthm}

\def\dbF{\hbox{\rm l\negthinspace F}}

\def\dbP{\hbox{\rm l\negthinspace P}}

\def\dbR{\hbox{\rm l\negthinspace R}}

\def\a{\alpha}
\def\b{\beta}

\def\d{\delta}

\def\cU{{\cal U}}

\def\sqr#1#2{{\vcenter{\vbox{\hrule height.#2pt
              \hbox{\vrule width.#2pt height#1pt \kern#1pt \vrule width.#2pt}
              \hrule height.#2pt}}}}

\def\dbR{{\mathop{\rm l\negthinspace R}}}

\def\3n{\negthinspace \negthinspace \negthinspace }
\def\2n{\negthinspace \negthinspace }
\def\1n{\negthinspace }

\def\dbF{{\mathop{\rm l\negthinspace F}}}

\def\dbP{{\mathop{\rm l\negthinspace P}}}
\def\dbR{{\mathop{\rm l\negthinspace R}}}
\def\={\buildrel \triangle \over =}

\def\a{\alpha}
\def\b{\beta}

\def\d{\delta}

\def\cU{{\cal U}}

\def\no{\noindent}

\def\ms{\medskip}

\def\max{\mathop{\rm max}}

\def\exp{\mathop{\rm exp}}
\def\sup{\mathop{\rm sup}}

\def\inf{\hbox{\rm inf$\,$}}

\def\as{\hbox{\rm a.s.{ }}}

\def\|{\Big |}
\def\({\Big (}
\def\){\Big )}
\def\[{\Big[}
\def\]{\Big]}
\def\be{\begin{equation}}
\def\bel{\begin{equation}\label}
\def\ee{\end{equation}}
\def\bt{\begin{theorem}}
\def\bcd{\begin{condition}}
\def\ecd{\end{condition}}
\def\et{\end{theorem}}
\def\bc{\begin{corollary}}
\def\ec{\end{corollary}}
\def\bde{\begin{definition}}
\def\ede{\end{definition}}
\def\bl{\begin{lemma}}
\def\el{\end{lemma}}
\def\bp{\begin{proposition}}
\def\ep{\end{proposition}}
\def\br{\begin{remark}}
\def\er{\end{remark}}
\def\ba{\begin{array}}
\def\ea{\end{array}}
\def\ed{\end{document}}

\def\square#1{\vbox{\hrule\hbox{\vrule height#1%
     \kern#1\vrule}\hrule}}
\def\rectangle#1#2{\vbox{\hrule\hbox{\vrule height#1%
     \kern#2\vrule}\hrule}}
\def\qed{\hfill \vrule height7pt width3pt depth0pt}

\font\tenbb=msbm10 \font\sevenbb=msbm7 \font\fivebb=msbm5

\newfam\bbfam
\scriptscriptfont\bbfam=\fivebb \textfont\bbfam=\tenbb
\scriptfont\bbfam=\sevenbb

\newtheorem{lemma}{Lemma}
\newtheorem{remark}{Remark}

\newtheorem{theorem}{Theorem}
\newtheorem{corollary}{Corollary}

\newtheorem{definition}{Definition}
\newtheorem{proposition}{Proposition}
\newtheorem{condition}{Assumption}

\begin{document}

\title{A Stochastic Maximum Principle for Risk-Sensitive Mean-Field Type Control}
\author{Boualem Djehiche$^{1}$, Hamidou Tembine$^{2}$ and Raul Tempone$^{2}$\\
$^{1}$ Department of Mathematics, KTH Royal Institute of Technology,\\ Stockholm, Sweden\\
 $^{2}$ SRI Uncertainty Quantification Center \\ in Computational Science and Engineering, CEMSE KAUST, KSA}
\maketitle
\begin{abstract}
In this paper we study mean-field type  control problems with risk-sensitive performance functionals. We establish a stochastic maximum principle (SMP) for optimal control of stochastic differential equations (SDEs) of mean-field type, in which the drift and the diffusion coefficients as well as the performance functional depend not only on the state and the control but also on the mean of the distribution of the state. Our result extends the risk-sensitive SMP (without mean-field coupling) of Lim and Zhou (2005), derived for feedback (or Markov) type optimal controls, to optimal control problems for non-Markovian  dynamics which may be time-inconsistent in the sense that the Bellman optimality principle does not hold. In our approach to the risk-sensitive SMP, the smoothness assumption on the value-function imposed in Lim and Zhou (2005) need not to be satisfied. For a general action space a Peng's type SMP is derived, specifying the necessary conditions for optimality. Two examples  are carried out to illustrate the proposed risk-sensitive mean-field type SMP under linear stochastic dynamics with exponential quadratic cost function. Explicit solutions are given for both mean-field free and mean-field models.
\end{abstract}

\begin{center}
\begin{minipage}{14.6cm}
\end{minipage}
\end{center}

\begin{center}
\begin{minipage}{14.6cm}
{\bf Abbreviated title.} Risk-sensitive control of SDEs of mean field type
\end{minipage}
\end{center}

\begin{center}
\begin{minipage}{14.6cm}
\end{minipage}
\end{center}

\section{Introduction}

Recently, there has been a renewed interest in optimal control problems for diffusions of mean-field type, where the performance functionals, drifts and diffusions coefficients depend not only on the state and the control but also on the probability distribution of state-control pair. Most formulations of mean-field type control in \cite{b1,b2,b3, hosking, li} have been of risk-neutral type where the performance functionals are the expected values of stage-additive payoff functions. Not all behavior, however, can be captured by risk-neutral mean-field type controls. One way of capturing risk-averse and risk-seeking behaviors is by exponentiating the performance functional before expectation (see \cite{Jac73}). 

A stochastic maximum principle (SMP) for the risk-sensitive optimal control problems for Markov diffusion processes with an exponential-of-integral performance functional was elegantly derived in \cite{lim} using the relationship between the SMP and the Dynamic Programming Principle (DPP) which expresses the first order adjoint process as the gradient of the value-function of the underlying control problem. This relationship holds only when the value-function is smooth (see Assumption (B4) in \cite{lim}). The approach of \cite{lim} was widely used and extended to jump processes in \cite{shi-1} and \cite{shi-2},  but still under this smoothness assumption. However, in many cases of interest, the value function is, in the best case, only continuous. Moreover, the relationship between the SMP and the DPP does not hold for non-Markovian dynamics and for mean-field type control problems where the Bellman optimality principle does not hold. This calls for the need to find a risk-sensitive SMP for these case. 

\ms The only paper that we are aware of and which deals with risk-sensitive optimal control in a mean field context is \cite{tembine2014}. Therein, the authors derive a verification theorem for a risk-sensitive mean-field game whose underlying dynamics is a Markov diffusion, using a matching argument between a system of Hamilton-Jacobi-Bellman (HJB) equations and the Fokker-Planck equation. This matching arguments freezes the mean-field coupling in the dynamics, which yields a standard risk-sensitive HJB equation for the value-function. The mean-field coupling is then retrieved through the Fokker-Planck equation satisfied by the marginal law of the optimal state. 

\ms Our contribution can be summarized as follows. We establish a stochastic maximum principle for a class of risk-sensitive mean-field type control problems where the distribution enters only through {\em the mean of state process}. This means that the drift, diffusion, running cost and terminal cost functions  depend on the state, the control and on the  mean of  state. Our work extends the results of \cite{lim} to risk-sensitive control problems for dynamics that are non-Markovian and of mean-field type. Our derivation of the SMP does not require any relationship between the first-order adjoint process and a value-function of an underlying control problem.  Using the SMP derived in \cite{hosking}, our approach can be easily extended to the case where the mean-field coupling is in terms of the mean of the state and the control processes.

\ms To the best to our knowledge, the risk-sensitive maximum principle for mean-field type controls has not been established in earlier work, is entirely new,  and is fundamentally different from the existing results in the risk-neutral mean-field case \cite{b1,b2,b3,hosking, li}.

\ms The paper is organized as follows. In Section \ref{sec:results}, we present the model and state the main result. In Section \ref{sec:intermediate}, we establish a risk sensitive SMP, based on the risk-neutral SMP by Buckdhan {\em et al.} \cite{b3}. In Section \ref{sec:mainresult}, we establish the risk-sensitive SMP. In section \ref{sec:example} we apply the risk-sensitive SMP to the linear-exponential- quadratic setup. Section \ref{sec:conclusion} concludes the paper. To streamline the presentation, we only consider the one-dimensional case. The extension to the multidimensional case is by now straightforward.

\section{Statement of the problem}\label{sec:results}

\hskip 1cm Let $T>0$ be a fixed time horizon and $(\Omega,{\cal{F}},\dbF, \dbP)$  be a given filtered
probability space on which a one-dimensional standard
Brownian motion $B=\{B_s\}_{s\geq0}$ is given, and the filtration $\dbF=\{{\mathcal{F}}_s,\ 0\leq s \leq T\}$ is the natural filtration of $B$ augmented by $\dbP-$null sets of ${\cal F}.$

 \ms We consider the stochastic control system:
\begin{equation}\label{SDEu}
\left\{\begin{array}{lll}
dx^u(t)=b(t,x^u(t),E[x^u(t)], u(t))dt+\sigma(t,x^u(t),E[x^u(t)], u(t))dB_t, \\ x^u(0)=x_0,
\end{array}\right.
\end{equation}
where
\begin{equation*}
b(t,x,y,u),\, \sigma(t,x,y,u): \,\,[0,T] \times \dbR\times\dbR\times U\longrightarrow \dbR,\ t\in[0, T],\ x\in \dbR,\ y\in \dbR,\ u \in U.
\end{equation*}

\ms\no An admissible control $u$ is an $\dbF$-adapted and square-integrable process with values in a non-empty subset $U$ of $\dbR^d$. We denote the set of all admissible controls by $\mathcal{U}$.

\ms\no Given $u\in\cU$ , equation (\ref{SDEu}) is an SDE with random coefficients. 

 \ms The risk-sensitive cost functional associated with (\ref{SDEu}) is given by
\begin{equation}\label{rs-cost}
J^{\theta}(u(\cdot))=Ee^{\theta\left[\int_0^Tf(t,x^u(t),E[x^u(t)], u(t))\,dt+ h(x^u(T),E[x^u(T)])\right]},
\end{equation}
where,  $\theta$ is the risk-sensitivity index, 
\begin{equation*}
f(t,x,y,u): \,\,[0,T] \times \dbR\times\dbR\times U\longrightarrow \dbR,\ \ h(x,y): \,\,\dbR\times\dbR\longrightarrow \dbR, \ t\in[0, T],\ x\in \dbR,\ y\in \dbR,\ u \in U.
\end{equation*}

\ms\no Any $\bar u(\cdot)\in {\cal U}$ satisfying
\be\label{rs-opt-u}
  J^{\theta}(\bar u(\cdot))=\inf_{u(\cdot)\in {\cal U}}J^{\theta}(u(\cdot))
\ee
is called a risk-sensitive optimal control. The corresponding state process, solution of (\ref{SDEu}), is denoted by $\bar x(\cdot):=x^{\bar u}(\cdot)$.

\ms\no  The optimal control problem we are concerned with is to characterize the pair $(\bar x,\bar u)$ solution of the  problem (\ref{rs-opt-u}).

\ms Let $\Psi_T=\int_0^T f(t, x(t),E[x(t)], u(t)) dt+h(x(T), E[x(T)])$. Then the risk sensitive loss functional is given by 
$$
\Psi_{\theta}:=\frac{1}{\theta}\log  E e^{\theta \left[\int_0^T f(t, x(t), E[x(t)],u(t))dt+h(x(T),  E[x(T)])\right]} = \frac{1}{\theta}\log \left[ E e^{\theta \Psi_T} \right].
$$
When the risk-sensitive index $\theta$ is small, the loss functional $\Psi_{\theta}$ can be expanded as
$$
E[ \Psi_T] +\frac{\theta}{2}\mbox{var}(\Psi_T)+O(\theta^2),
$$
where,  $\mbox{var}(\Psi_T)$ denotes the variance of  $\Psi_T$. If $\theta<0$ , the variance of $\Psi_T$, as a measure of risk,  improves the performance $\Psi_{\theta}$, in which case the optimizer is called {\it risk seeker}. But,  when $\theta>0$,  the variance of $\Psi_T$ worsens the performance $\Psi_{\theta}$,   in which case the optimizer is called {\it risk averse}.
The risk-neutral loss functional $E[ \Psi_{T}]$ can be seen as a limit of risk-sensitive functional $ \Psi_{\theta}$ when $\theta\rightarrow 0$.

\ms\no Note that the presence of the expectations  $ E[x(T)]$  in the loss function  $\Psi_{T}$  may cause  time-inconsistency, in which case the Bellman's Principle is no longer valid and this motivates the use of the stochastic maximum (SMP) approach instead of  trying extensions of the dynamic programming principle (DPP).

 \ms For convenience, we will use the following notation throughout the paper. For $\phi\in \{b, \sigma, f, h\}$, etc, respectively, we define
\be\label{notation}\left\{\begin{array}{llll}
\delta\phi(t)=\phi(t,\bar x(t),E[\bar x(t)],u(t))-\phi(t,\bar x(t),E[\bar x(t)],\bar u(t));\\
\phi_x(t)=\frac{\partial\phi}{\partial x}(t,\bar x(t),E[\bar x(t)],\bar u(t)),\,\,\,\phi_{xx}(t)=\frac{\partial^ 2\phi}{\partial x^ 2}(t,\bar x(t),E[\bar x(t)],\bar u(t));\\
\phi_y(t)=\frac{\partial\phi}{\partial y}(t,\bar x(t),E[\bar x(t)],\bar u(t)).
\end{array}\right.
\ee
where $u$ is an admissible control from ${\cal U}$.

\ms\noindent We define the risk-neutral Hamiltonian associated with random variables $X\in L^1(\Omega, {\cal F},\dbP)$ as follows.  for $(p, q)\in \dbR\times \dbR$
\begin{equation}\label{hamiltonian}
H(t,X,u,p,q):=b(t,X,E\left[X\right],u)p+\sigma(t,X,E\left[X\right],u)q-f(t,X,E\left[X\right],u),
\end{equation}

\no We also introduce the risk-sensitive Hamiltonian: for $\theta \in \dbR$ and $(p,q, \ell)\in \dbR\times \dbR\times \dbR$,
\begin{equation}\label{rs-hamiltonian}
H^{\theta}(t, X, u, p, q,\ell):= b(t, X,E\left[X\right], u) p+\sigma(t,X,E\left[X\right], u) (q +\theta\ell p)-f (t, X,E\left[X\right], u).
\end{equation}
We have $H=H^0$.

\ms\no  Moreover, we  denote
\be\label{dH}\begin{array}{lll}
\d H(t):=p(t)\d b(t)+q(t)\d \sigma(t)-\d f(t), \quad \d H^{\theta}(t):=p(t)\d b(t)+(q +\theta\ell p)\d \sigma(t)-\d f(t),\\
H_k(t):=b_k(t)p+\sigma_k(t)q-f_k(t),\quad H_k^{\theta}(t):=p(t)b_k(t)+(q +\theta\ell p) \sigma_k(t)- f_k(t),
\end{array}
\ee
for $k=x,y, xx.$

\ms We will make the following assumption in this paper.
\bcd\label{cond1} The functions $b, \sigma, f, h$ are twice continuously differentiable with respect to $(x,y)$.
Moreover, $b, \sigma, f,h$ and all their derivatives up to second order with respect to $(x,y)$ are continuous in $(x,y, u)$, and bounded.
\ecd

\no  Under these assumptions, for each $u\in {\cal U}$, the SDE (\ref{SDEu}) admits a unique strong solution $x^u$ (see e.g. \cite{b3, Buckdahn1}).

\ms\no 
We introduce the adjoint equations involved in the risk-sensitive SMP for our control problem.

\ms
The first order adjoint equation is the following  backward SDE of mean-field type:
\be\label{rs-firstAD-1}\left\{\begin{array}{lll}
d \bar{p}(t)= -\left\{H^{\theta}_x(t)+\frac{1}{v^{\theta}(t)}E[v^{\theta}(t) H^{\theta}_y(t)]\right\} dt + \bar{q}(t)(-\theta\ell(t)dt+ d B_t),
\\
dv^{\theta}(t)=\theta\ell(t)v^{\theta}(t)dB_t,
\\
v^{\theta}(T)=\phi^{\theta}(T),\\
\bar{p}(T)= -h_x(T)-\frac{1}{\phi^{\theta}(T)}E[\phi^{\theta}(T) h_{y}(T)].
\end{array}
\right.
\ee
where,
\be\label{phi}
\phi^{\theta}_T:=e^{\theta [h(\bar x(T), E[\bar x(T)])+ \int_0^Tf(t,\bar x(t), E[\bar x(t)], \bar u(t)) dt]}.
\ee
In view of (\cite{Buckdahn1}, Theorem 3.1.), under Assumption 1, (\ref{rs-firstAD-1}) admits a unique $\dbF$-adapted solution $(\bar p,\bar q,v^{\theta},\ell)$ such that
\be\label{rs-first-bounds}
E\left[\sup_{t\in[0,T]}|\bar p(t)|^2+\sup_{t\in[0,T]}|v^{\theta}(t)|^2+\int_0^T \left(|\bar q(t)|^2+|\ell(t)|^2\right) dt\right]<\infty.
\ee

\ms\no The second order adjoint equation is the following backward SDE: 

\be\label{rs-secondAD-1}\left\{ \begin{array}{lll}
d\bar{P}(t)=-\left\{ \left(2b_x(t)+\sigma_x^2(t)+2\theta\ell(t)\sigma_x(t)\right)\bar{P}(t)+2\sigma_x(t)\bar{Q}(t)\right. \\ \qquad\quad  \left. - \theta(\bar{q}(t)+\sigma_x(t)\bar{p}(t))^2+H^{\theta}_{xx}(t)\right\}dt+\bar{Q}(t)(-\theta\ell(t)dt+dB_t),\\ 
\bar{P}(T)=-h_{xx}(T).
\end{array}
\right.
\ee
This is a standard linear backward SDE, whose unique $\dbF$-adapted solution $(\bar P,\bar Q)$ satisfies
\be\label{rs-second-bounds}
E\left[\sup_{t\in[0,T]}|\bar P(t)|^2+\int_0^T |\bar Q(t)|^2 dt\right]<\infty.
\ee 
 
\no The following theorem is the main result of the paper.

\begin{theorem}\label{main result}$(${\bf Risk-sensitive maximum principle}$)$
Let Assumption \ref{cond1} hold. If  $(\bar x(\cdot),\bar  u(\cdot))$  is an optimal solution of the risk-sensitive control problem (\ref{SDEu})-(\ref{rs-cost}), then there are three pairs of $\dbF$-adapted processes $(v^{\theta}, \ell)$, $(\bar p,\bar{q})$ and $(\bar{P}, \bar{Q})$ that satisfy (\ref{rs-firstAD-1})-(\ref{rs-first-bounds}) and (\ref{rs-secondAD-1})-(\ref{rs-second-bounds}) respectively, such that
\be\label{rs-VI-1}
\d H^{\theta}(t)+\frac{1}{2}\left(\bar{P}(t)-\theta \bar{p}^2(t)\right)\left(\d\sigma(t)\right)^2\le 0,
\ee
for all $u\in U$, almost every $t\in [0,T]$ and $\mathbb{P}-$almost surely.

\medskip
\no In particular, if $\sigma(t,x,u):=\sigma(t,x)$ is independent of $u$ then 

$$
H^{\theta}(t, \bar x(t),\bar u(t), \bar{p}(t), \bar{q}(t),\ell(t))=\max_{u} H^{\theta}(t,\bar x(t),  u, \bar{p}(t), \bar{q}(t),\ell(t)).
$$
\end{theorem}

\begin{remark}  Theorem \ref{main result} reduces to Theorem 3.1  of Lim and Zhou \cite{lim}, if the model is mean-field free  i.e. for which $
\sigma_{y}=h_{y}=b_{y}=f_{y}=0$, and when  $\ell(t):=-\bar p(t)\sigma(t, \bar x(t), \bar u(t))$, in which case the generic martingale $v^{\theta}$ becomes the smooth value-function of a  Markovian or feedback control dynamics, whose gradient is  the adjoint process.
\end{remark}

\begin{remark}
The main results of the paper are built on the SMP for the risk neutral case derived in \cite{b3}, where the strong condition \ref{cond1} on the involved coefficients is imposed to get less technical proofs.  These conditions  can be considerably weakened using techniques that are by now well established in the optimal control literature (see e.g.  \cite{Mezerdi,lim}). 
\end{remark}

\section{Proof of the risk-sensitive stochastic maximum principle}

\no The proof of Theorem \ref{main result} is displayed in the next subsections.

\subsection{An intermediate SMP for mean-field type control}\label{sec:intermediate}
 In this subsection we first reformulate the risk-sensitive control problem (\ref{SDEu})-(\ref{rs-opt-u}) in terms of an augmented state process and terminal payoff  problem. An
intermediate stochastic maximum principle is then obtained by applying the SMP of (\cite{b3}, Theorem 2.1.) for loss functionals without running cost. Then,
we transform the intermediate first- and second-order adjoint processes to a more simpler form.  The  mean-field type control  problem  (\ref{rs-opt-u}) under the dynamics (\ref{SDEu}) is equivalent to 
\begin{eqnarray}
\label{smp2}
\left\{
\begin{array}{lll} \inf_{u(\cdot)\in \cU}  E e^{\theta \left[h(x(T), E[x(T)])+ \xi(T)\right]},
\\
\displaystyle{\mbox{ subject to }\ }\\
dx(t)=b(t,x(t),E[x(t)], u(t))dt+\sigma(t,x(t),E[x(t)], u(t))dB_t,\\
d\xi(t)=f(t,x(t), E[x(t)], u(t)) dt,\\
x(0)=x_{0}, \quad \xi(0)=0. \\
 \end{array}
\right.
\end{eqnarray}
Recall that
 \begin{equation*}\label{phi}
\phi^{\theta}_T:=e^{\theta [h(\bar x(T), E[\bar x(T)])+ \int_0^Tf(t,\bar x(t), E[\bar x(t)], \bar u(t)) dt]}.
\end{equation*}
Under Assumption \ref{cond1}, we may apply the  SMP for risk-neutral mean-field type control from  (\cite{b3}, Theorem 2.1) to the augmented state dynamics $(x,\xi)$ to derive the first order adjoint equation

\begin{equation}\label{rn-firstAD} \left\{ \begin{array}{lll}  
 d\vec{p}(t)= -\left\{ \left( \begin{array}{cc}
b_x(t) & 0  \\
f_x(t) & 0   
\end{array} \right)^{\prime} \vec{p}(t)+ \left( \begin{array}{cc}
\sigma_x(t) & 0  \\
0 & 0   
\end{array} \right)^{\prime} \vec{q}
+E \left[\left(\begin{array}{cc}
b_{y}(t) & 0  \\
f_{y}(t) & 0   
\end{array} \right)^{\prime}\vec{p}
+\left( \begin{array}{cc}
\sigma_{y}(t) & 0  \\
0 & 0   
\end{array} \right)^{\prime} \vec{q}(t)\right]\right\}dt\\ \\ \qquad\qquad + \vec{q}(t) dB_t,\\ \\
\vec{p}(T)=-\theta \phi^{\theta}_T\left(\begin{array}{cc}
h_{x}(T)   \\
1    \end{array} \right) -   \theta E\left[ \phi^{\theta}_T\left(\begin{array}{cc}
h_{y}(T)   \\ 0
\end{array} \right)\right],
\end{array}
\right.
\end{equation} 
with
\be\label{rn-pq-ineq}
E\left[\sup_{0\le t\le T}|\vec{p}(t))|^2+\int_0^T|\vec{q}(t)|^2 dt \right ]<\infty.
\ee

Let ${\tilde H}^{\theta}$  be the Hamiltonian associated with the optimal state dynamics $\bar x(\cdot)$ and the pair of adjoint processes $( \vec{p}(t), \vec{q}(t))$:
\be\label{inter-ham}
 {\tilde H}^{\theta}(t, \bar x(t), u, \vec{p}(t), \vec{q}(t)):=  \left(\begin{array}{cc}
b(t, \bar x(t), E[\bar x(t)], u)   \\
f(t,  \bar x(t), E[\bar x(t)], u)   \end{array} \right)\cdot \vec{p}(t)+\left(\begin{array}{cc}
\sigma(t,  \bar x(t), E[\bar x(t)], u) \\ 0  \end{array}\right)\cdot\vec{q}(t),
\ee
where, $(\,\cdot\,)$ denotes the usual scalar product in $\dbR^2$.  The dependence of the Hamiltonian on $\theta$ stems from the dependence of the adjoint processes $(\vec{p},\vec{q})$ of $\theta$ through the end-condition in (\ref{rn-firstAD}). 

\ms\no The second order adjoint equation is 
\begin{equation}\label{rn-secondAD} 
\left \{ \begin{array}{lll}
dP(t)= -\left\{ \left( \begin{array}{cc}
b_x(t) & 0  \\
f_x(t) & 0   
\end{array} \right) P(t) +P(t)\left( \begin{array}{cc}
b_x(t) & 0  \\
f_x(t) & 0   
\end{array} \right)^{\prime}
+\left( \begin{array}{cc}
\sigma_x(t) & 0  \\ 0 & 0   
\end{array} \right) P(t) \left( \begin{array}{cc}
\sigma_x(t) & 0  \\
0 & 0   
\end{array} \right)^{\prime}  \right. \\ \\  \quad \qquad \left. +
\left( \begin{array}{cc}
\sigma_x(t) & 0  \\
0 & 0   
\end{array} \right)Q(t)  +  Q(t)\left( \begin{array}{cc}
\sigma_x(t) & 0  \\
0 & 0   
\end{array} \right)^{\prime} + \left( \begin{array}{cc}
{\tilde H}_{xx}^{\theta}(t) & 0  \\ 
0 & 0   
\end{array} \right) \right \}dt+ Q(t) dB_t, \\ \\ 
P(T)=-\theta \phi_T
\left( \begin{array}{cc}
\theta h_x^2(T)+ h_{xx}(T) & \theta h_x(T)  \\
\theta h_x(T) & \theta
\end{array} \right).
\end{array}
\right.
\end{equation} 
\be\label{rn-PQ-ineq}
E\left[\sup_{0\le t\le T}||P(t))||^2+\int_0^T||Q(t)||^2dt \right ]<\infty,
\ee
where, $||\cdot||$ denotes the norm of the coresponding  matrices.
 
\ms\no We have the following 
\begin{proposition}\label{rn-MSP} 
Let Assumption \ref{cond1} hold. If $(\bar x,\bar\xi, \bar  u)$ is an optimal solution of the risk-neutral control problem  (\ref{smp2}),  then there are two pairs of $\dbF$-adapted processes $(\vec{p},\vec{q})$ and $(P, Q)$ that satisfy (\ref{rn-firstAD})-(\ref{rn-pq-ineq}) and (\ref{rn-secondAD})-(\ref{rn-PQ-ineq}) respectively, such that
\begin{equation}\label{rn-VI}\begin{array}{lll}
\d{\tilde H}^{\theta}(t)+\frac{1}{2}\left( \begin{array}{cc}
\d\sigma(t)  \\ 0   
\end{array} \right)^{\prime}P(t)\left( \begin{array}{cc}
\d\sigma(t) \\ 0   
\end{array} \right) \leq 0, 
\end{array}
\end{equation} for all $u\in U,$ almost every $t$ and $\mathbb{P}-$almost surely,

\medskip\noindent 
where,
$$
\d{\tilde H}^{\theta}(t):={\tilde H}^{\theta}(t, \bar x(t), u,\vec{p}(t), \vec{q}(t))-{\tilde H}^{\theta}(t, \bar x(t),  \bar u(t), \vec{p}(t) ,\vec{q}(t)).
$$
\end{proposition}

\subsection{Transformation of the first order adjoint process}\label{sec:mainresult}
Although the result of Proposition \ref{rn-MSP} is a good SMP for the risk-sensitive mean-field type control, the fact that augmenting the state process with the second component $\xi$  yields a system of two  adjoint equations that appears  complicated to solve in concrete situations. In the mean-field free case,  Lim and Zhou (\cite{lim}) elegantly  solve this problem by  suggesting a transformation of the adjoint processes $(\vec{p},\vec{q})$ in  such a way to get rid of the second component $(p_2,q_2)$ in  (\ref{rn-firstAD}) and express the SMP in terms of only one adjoint process that we denote $(\tilde p_1,\tilde q_1)$ and which solves a backward SDE whose driver is quadratic in $\tilde{p_1}$, which reminds of the risk-sensitive Hamilton-Jacobi-Bellman equation (see \cite{tembine2014} and the references therein). The suggested transform uses a relationship between the SMP  and  the DPP  (valid only for Markovian or feedback controls and in situations where the Bellman Principle is valid) which expresses the adjoint process $\vec{p}$ as the gradient of the value-function associated with the control problem (\ref{smp2}), provided that the value-function is smooth (see Assumption (B4) in \cite{lim}), a condition that is  often hard to verify in concrete situations. The value-function is in general {\it not} smooth. Furthermore, the approach developed in \cite{lim} cannot be extended general situations, such as non-Markovian dynamics and in mean-field type control problems, where the Bellman Principle does not hold.

\ms\no A closer look at the method of Lim and Zhou (\cite{lim}), suggests in fact that it is enough to use  a  generic square-integrable martingale to transform  the pair $(p_2, q_2)$  into the adjoint process $(\tilde{p}_2, 0)$, where the process $\tilde{p}_2$ is still a square-integrable martingale, which would mean that $\tilde{p}_2(t)=\tilde{p}_2(T)$ and is equal to the  constant $E[\tilde{p}_2(T)]$. But, this generic martingale need not be related to the adjoint process $\vec{p}$ as in (\cite{lim}). Instead,  it will be part of the adjoint equation associated with the risk-sensitive SMP (see Theorem \ref{main result}, above).

\medskip\noindent  Roughly,  noting that 
$dp_2(t) = q_2(t) dB_t$ and $ p_2(T)=- \theta \phi^{\theta}_T,$  the explicit solution of this backward SDE is  
\be\label{p2}
p_2(t)=- \theta E[\phi^{\theta}_T \ | \ {\cal F}_t ]=-\theta v^{\theta}(t) ,
\ee
where, 
\be\label{mg} 
v^{\theta}(t):= E[\phi^{\theta}_T \ | \ {\cal F}_t ],  \qquad 0\le t\le T. 
\ee
In view of (\ref{p2}), it would be  natural to choose a transformation of $(\vec{p},\vec{q})$ into an adjoint process $(\tilde{p}, \tilde{q})$ ,
where,
$$
\tilde p(t):=\left( \begin{array}{cc}
\tilde p_1 (t)\\ \tilde p_2(t)   
\end{array} \right), \,\,\,   \tilde q(t):=\left( \begin{array}{cc}
\tilde q_1 (t)\\ \tilde q_2(t)   
\end{array} \right), \quad 0\le t\le T, 
$$
 such that
\be\label{t-p2}
\tilde p_2(t)=\frac{p_2(t)}{ \theta v^{\theta}(t)}=-1, \qquad 0\le t\le T,
\ee
which would  imply that, for almost  every $0\le t\le T$, $\tilde q_2(t)=0, \quad \dbP-\as$ 

\medskip We consider the following transform  
\be\label{transform}
\tilde p(t)=\left( \begin{array}{cc}
\tilde p_1 (t)\\ \tilde p_2(t)   
\end{array} \right)  :=\frac{1}{\theta v^{\theta}(t)} \vec{p}(t),\qquad  0\le t\le T.
\ee
In view of (\ref{rn-firstAD}) and (\ref{v}), we have
\be\label{t-p-T}
\tilde{p}(T)= - \left( \begin{array}{cc}
h_x(T)+\frac{1}{v^{\theta}(T)}E[v^{\theta}(T)h_y(T)] \\
1 
\end{array} \right).
\ee
The following properties  of the generic martingale $v^{\theta}$ are essential in order to investigate the properties of these new processes $(\tilde p,\tilde q)$.

\ms\no First, we note that since, by Assumption \ref{cond1}, $f$ and $h$ are bounded by some constant $C>0$, we have 
\be\label{v-bound}
0< e^{-(1 + T)C\theta} \leq \phi^{\theta}_T \leq  e^{(1 + T)C\theta}.
\ee
Therefore, $v^{\theta}$  is a uniformly bounded $\dbF$-martingale satisfying
\be\label{b-mg}
0< e^{-(1 + T)C\theta} \leq v^{\theta}(t) \leq e^{(1 + T)C\theta}, \,\qquad  0\le t\le T.
\ee
 Furthermore, the martingale $v^{\theta}$ enjoys the following useful {\it  logarithmic transform}  established in (\cite{Karoui-Ham}, Proposition 3.1): 
\be\label{mart-rep}
v^{\theta}(t)=\exp\left(\theta Y_t+\theta \int_0^t f(s,\bar x(s), E[\bar{x}(s)],\bar u(s))ds\right), \quad 0\le t\le T, 
\ee 
and 
\be\label{Y-v}
v^{\theta}(0)=E[\phi^{\theta}_T]=\exp(\theta Y_0),
\ee
where,  in view of (\ref{b-mg}) and the boundedness of $f$, 
\be\label{b-Y}
\sup_{0\le t\le T }|Y_t|\le C_T,
\ee
where, $C_T$ is  a positive constant that depends only on $T$ and the bounds of $f$ and $h$. Moreover, the process $Y$ is the first component of the $\dbF$-adapted pair of processes $(Y,\ell)$  which is the unique solution to the following quadratic BSDE:
\be\label{Q-BSDE}
\left\{ \begin{array}{lll}
dY_t=-\{f(t, \bar x(t), E[\bar{x}(t)] ,\bar u(t))+\frac{\theta}{2}|\ell(t)|^2\}dt+\ell(t)dB_t,\\ 
Y_T=h(\bar x_T, E[\bar x_T)]),\\
\end{array}
\right.
\ee
where,
\be\label{Y-ell}
E\left[\int_0^T |\ell(t)|^2dt\right]<\infty.
\ee

\medskip\noindent 
In particular, $v^{\theta}$ solves the following linear backward SDE
\be\label{v}
dv^{\theta}(t)=\theta\ell(t)v^{\theta}(t)dB_t,\quad  v^{\theta}(T)=\phi^{\theta}_T.
\ee
Hence, 
\be\label{v-L}
\frac{v^{\theta}(t)}{v^{\theta}(0)}=\exp{\left(\int_0^t \theta\ell (s)dB_s - \frac{\theta^2}{2}\int_0^t |\ell (s)|^2 ds\right)}:=L^{\theta}_t, \quad 0\le t\le T.
\ee
is a uniformly bounded $\dbF$-martingale.

\ms\no We wish to identify the processes $\tilde\a$ and $\tilde q$ such that
\be\label{t-p}
d\tilde p(t)=-\tilde\a(t)dt+\tilde q(t) dB_t.
\ee
We may  apply It\^o's formula to  the process $\vec{p}(t)=\theta v^{\theta}\tilde{p}(t)$, using the expression of $ v^{\theta}$ in (\ref{v}), to obtain  
$$
d\vec{p}(t)= \theta v^{\theta}(t) d\tilde{p}(t) + \theta^2 \ell(t)v^{\theta}(t) \tilde{q}(t)dt +
  \theta^2 \ell(t)v^{\theta}(t)\tilde{p}(t) dB_t.
$$
Thus,
$$ 
d\tilde{p}(t)=\frac{d\vec{p}(t)}{\theta v^{\theta}(t)}   -  \theta\ell(t)\tilde{q}(t) dt -
\tilde{p}(t) \theta\ell(t) dB_t.
$$
Substituting the expression of $d\vec{p}$ in (\ref{rn-firstAD}) and identifying the coefficients we get the diffusion term
\begin{eqnarray}\label{t-q}
\tilde{q}(t)= \frac{1}{\theta  v^{\theta}(t)}\vec{q}(t) -\theta\ell(t)\tilde{p}(t), \quad 0\le t\le T,
\end{eqnarray}
and the drift term of the process $\tilde{p}$ 
\begin{equation*}\begin{array}{lll}
\tilde{\alpha}(t)=  
 \left( \begin{array}{cc}
b_x(t) & 0  \\
f_x(t) & 0   
\end{array} \right)^{\prime} \frac{1}{\theta v^{\theta}(t)}\vec{p}(t)+ \left( \begin{array}{cc}
\sigma_x (t)& 0  \\
0 & 0   
\end{array} \right)^{\prime} \frac{1}{\theta v^{\theta}(t)}\vec{q}(t) \\ \\ \qquad 
+\frac{1}{\theta v^{\theta}(t)}E\left[\left(\begin{array}{cc}
b_{y}(t) & 0  \\
f_y(t) & 0   
\end{array} \right)^{\prime}\vec{p}(t)+ \left( \begin{array}{cc}
\sigma_{y}(t) & 0  \\
0 & 0   
\end{array} \right)^{\prime} \vec{q}\right]
+\theta\ell(t)\tilde{q}(t).
\end{array}
\end{equation*}
Now using the relations 
$$
 \vec{p}(t)= \theta v^{\theta}(t) \tilde{p}(t),\quad 
\vec{q}(t)= \theta v^{\theta}(t) \tilde{q}(t) +\theta \tilde{p}(t) \theta\ell(t)v^{\theta}(t),
$$
we finally obtain
\be \label{t-alpha}\begin{array}{lll}
\tilde{\alpha}(t)=
 \left( \begin{array}{cc}
b_x(t) & 0  \\ 
f_x(t) & 0   
\end{array} \right)^{\prime} \tilde{p}(t) + \left( \begin{array}{cc}
\sigma_x (t)& 0  \\
0 & 0  \end{array} \right)^{\prime} \left\{\tilde{q}(t) +  \theta\ell(t)\tilde{p}(t) \right\} + \theta\ell(t)\tilde{q}(t) \\ \\ \qquad
+ \frac{1}{v^{\theta}(t)}E\left [v^{\theta}(t)\left(\begin{array}{cc}
 b_{y}(t) & 0  \\ 
f_{y}(t) & 0   
\end{array} \right)^{\prime} \tilde{p}(t)
 + v^{\theta}(t) \left( \begin{array}{cc}
\sigma_{y}(t) & 0  \\
0 & 0   
\end{array} \right)^{\prime}  \left\{\tilde{q}(t)+  \theta\ell(t)\tilde{p}(t)  \right\} \right].
\end{array}
\ee
It is easily verified that
$$
d\tilde{p}_2(t)= \tilde{q}_2(t)[-\theta\ell(t) dt+  dB_t],\quad   \tilde{p}_2(T)=-1.
$$
In view of (\ref{v-L}),  we may use Girsanov Theorem to claim that
$$
d\tilde{p}_2(t)=\tilde{q}_2 (t)dB_t^{\theta},  \quad  \mathbb{P}^{\theta}-\as \qquad \tilde p_2(T)=-1,
$$
where,
$$
B_t^{\theta}:=B_t-\int_0^t \theta \ell(s)ds
$$
is a $\mathbb{P}^{\theta}$-Brownian motion, where,
$$ 
\frac{d\mathbb{P}^{\theta}}{d\mathbb{P}} \Large|_{{\cal F}_t} :=L^{\theta}_t=\exp{\left(\int_0^t \theta\ell (s)dB_s - \frac{\theta^2}{2}\int_0^t |\ell (s)|^2 ds\right)},\quad 0\le t\le T.
$$
In view of (\ref{v-L}) and (\ref{v-bound}), the probability measures $\dbP$ and $\dbP^{\theta}$ are in fact equivalent.
\no 
Hence, noting that $\tilde p_2(t):=[\theta v^{\theta}(t)]^{-1}p_2(t)$ is square-integrable, we get that  $\tilde p_2(t)=E^{\mathbb{P}^{\theta}}[\tilde p_2(T)| {\cal F}_t]=-1$.  Thus, its quadratic variation $\int_0^T|\tilde q_2(t)|^2dt=0$. This 
 implies that,  for almost every $0\le t\le T$, $\tilde q_2(t)=0,\,\,  \mathbb{P}^{\theta}\,\,\mbox{and}\,\, \mathbb{P}-\as$

\medskip\noindent  Therefore, the first component of $\tilde\a$ given by (\ref{t-alpha}) reads
\begin{equation} \begin{array}{lll}
\tilde{\alpha}_1(t)=b_x(t)\tilde{p}_1(t) +  \sigma_x(t) \left(\tilde{q}_1(t)+\theta\ell(t)\tilde{p}_1(t)\right) -f_x(t)+ \tilde{q}(t)\theta\ell(t) \\ \\ \qquad\quad +\frac{1}{v^{\theta}(t)}E\left[v^{\theta}(t)\left(b_y(t)\tilde{p}_1(t) +  \sigma_y(t) \left(\tilde{q}_1(t)+\theta\ell(t)\tilde{p}_1(t)\right) -f_y(t) \right) \right].
\end{array}
\end{equation}
and the main risk-sensitive first order adjoint equation for $(\tilde p_1,\tilde q_1)$ and $(v^{\theta},\ell)$ becomes
\be\label{t-p-1}\left\{\begin{array}{lll}
d \tilde{p}_1= -\tilde{\alpha}_1(t) dt + \tilde{q}_1 d B_t, \\
dv^{\theta}(t)=\theta\ell(t)v^{\theta}(t)dB_t,\\
v^{\theta}(T)=\phi^{\theta}(T),\\
\tilde{p}_1(T)= -h_x(T)-\frac{1}{\phi^{\theta}(T)}E[\phi^{\theta}(T) h_{y}(T)].
\end{array}
\right.
\ee
The solution  of this system of backward SDEs is unique.

\medskip
\subsection{Transformation of the Hamiltonian} 

In view of (\ref{transform}) and (\ref{t-q}), the Hamiltonian  $\tilde H^{\theta}$, associated with (\ref{smp2}),  given by

$$ 
 {\tilde H}^{\theta}(t, \bar x(t), u, \vec{p}(t), \vec{q}(t))=  \left(\begin{array}{cc}
b(t,  \bar x(t), E[\bar x(t)], u)   \\
f(t,  \bar x(t), E[\bar x(t)], u)   \end{array} \right)\cdot \vec{p}(t)+\left(\begin{array}{cc}
\sigma(t,  \bar x(t), E[\bar x(t)], u) \\ 0  \end{array}\right)\cdot\vec{q}(t),
$$ 
satisfies
\be\label{H-tilde-H}
{\tilde H}^{\theta}(t, \bar x(t), u, \vec{ p}(t),\vec{q}(t))=[\theta v^{\theta}(t)] 
H^{\theta}(t, \bar x(t),  u, \tilde{p}_1(t), \tilde{q}_1(t),\ell(t))
\ee
where, $ H^{\theta}$ is  the risk-sensitive Hamiltonian given by (\ref{rs-hamiltonian}). 

\no Using the notation (\ref{notation}), we have the following relation between the drift term $\tilde\a_1$ in (\ref{t-alpha}) and the gradient of the risk-sensitive Hamiltonian $H^{\theta}$:
\be\label{alpha-tilde-H}
\tilde\a_1(t)= H^{\theta}_x(t)+\frac{1}{v^{\theta}(t)}E[v^{\theta}(t) H^{\theta}_y(t)]+\theta\ell(t)\tilde{q}_1(t).
\ee  
Hence, risk-sensitive first order adjoint equation (\ref{t-p-1}) becomes
\be\label{rs-firstAD}\left\{\begin{array}{lll}
d \tilde{p}_1= -\left\{H^{\theta}_x(t)+\frac{1}{v^{\theta}(t)}E[v^{\theta}(t) H^{\theta}_y(t)]\right\} dt + \tilde{q}_1(-\theta\ell(t)dt+ d B_t),
\\
dv^{\theta}(t)=\theta\ell(t)v^{\theta}(t)dB_t,
\\
v^{\theta}(T)=\phi^{\theta}(T),\\
\tilde{p}_1(T)= -h_x(T)-\frac{1}{\phi^{\theta}(T)}E[\phi^{\theta}(T) h_{y}(T)].
\end{array}
\right.
\ee

\subsection{Transformation of the second order adjoint process}
For the second order adjoint equation, we apply  the same type of transformations suggested  in \cite{lim} and let  
\be\label{t-P}
\tilde{P}(t):=\frac{P(t)}{\theta v^{\theta}(t)}+\theta \tilde{p}(t) \tilde{p}^{\prime}(t),\quad 0\le t\le T, 
\ee
and
\be\label{t-Q}
\tilde{Q}(t):=\frac{Q(t)}{\theta v^{\theta}(t)}+\theta \tilde q(t) \tilde{p}^{\prime}(t)+\theta \tilde p(t) \tilde{q}^{\prime}(t)
-\theta\ell(t)\left( \tilde{P}(t)-\theta \tilde p(t) \tilde{p}^{\prime}(t)\right), \quad 0\le t\le T. 
\ee
In view of (\ref{rn-secondAD}) satisfied by $(P,Q)$ , easy (but lengthy) calculations similar to \cite{lim}  yield that  
\be\label{t-PQ-1}
\tilde{P}(t)= \left( \begin{array}{cc}
\tilde{P}_1 (t)& 0  \\
0 & 0  \end{array} \right) ,\quad \tilde{Q}= \left( \begin{array}{cc}
\tilde{Q}_1 (t)& 0  \\
0 & 0  \end{array} \right), \quad 0\le t\le T,
\ee
where, $(\tilde{P}_1, \tilde{Q}_1)$ is a pair of one-dimensional processes that uniquely solve the risk-sensitive second order adjoint equation:
\be\label{rs-secondAD}\left\{ \begin{array}{lll}
d\tilde{P}_1=-\left\{ \left(2b_x(t)+\sigma_x^2(t)+2\theta\ell(t)\sigma_x(t)\right)\tilde{P}_1(t)+2\sigma_x(t)\tilde{Q}_1(t)\right. \\ \\ \qquad\quad  \left. - \theta(\tilde{q}_1(t)+\sigma_x(t)\tilde{p}_1(t))^2+ H^{\theta}_{xx}(t)\right\}dt+\tilde{Q}_1(t)(-\theta\ell(t)dt+dB_t),\\ \\
\tilde{P}_1(T)=-h_{xx}(T).
\end{array}
\right.
\ee

\subsection{ Risk-sensitive stochastic maximum principle}

To arrive at a risk-sensitive SMP expressed in terms of the adjoint processes $(\tilde{p}_1,\tilde{q}_1), (v^{\theta}, \ell)$ and $(\tilde{P}_1,\tilde{Q}_1,)$, which solve (\ref{rs-firstAD}) and (\ref{rs-secondAD}) respectively, we note that in view of (\ref{t-P}), (\ref{t-Q}) and (\ref{t-PQ-1}),  the second term in  the variational inequality (\ref{rn-VI}) satisfies 
\begin{equation*}\begin{array}{lll}
\left( \begin{array}{c}
\d\sigma(t)   \\
0    
\end{array} \right)^{\prime} P(t)\left( \begin{array}{c}
\d\sigma(t)  \\
0    
\end{array} \right) =[\theta v^{\theta}(t)] \left( \begin{array}{c}
\d\sigma(t)   \\
0    
\end{array} \right)^{\prime} \left[ \tilde{P}(t) - \theta\tilde{p}(t)\tilde{p}^{\prime}(t)\right] \left( \begin{array}{c}
\d\sigma(t)  \\
0    
\end{array} \right) \\ \qquad\qquad\qquad \qquad\qquad\qquad\quad = {\color{black}\theta v^{\theta}(t)}\left(\tilde{P}_1(t)-\theta \tilde{p}^2_1(t)\right)\left(\d\sigma(t)\right)^2.
\end{array}
\end{equation*}
Combining this relation with  (\ref{H-tilde-H}) we obtain
$$
\d {\tilde H}^{\theta}(t)+\frac{1}{2}\left( \begin{array}{cc}
\d\sigma(t)  \\ 0   
\end{array} \right)^{\prime}P(t)\left( \begin{array}{cc}
\d\sigma(t) \\ 0   
\end{array} \right)=\theta v^{\theta}(t)\left[\d H^{\theta}(t)+\frac{1}{2}\left(\tilde{P}_1(t)-\theta \tilde{p}^2_1(t)\right)\left(\d\sigma(t)\right)^2\right].
$$
Hence, since $v^{\theta}>0$, the variational inequality  (\ref{rn-VI}) translates into
$$
\d H^{\theta}(t)+\frac{1}{2}\left(\tilde{P}_1(t)-\theta \tilde{p}^2_1(t)\right)\left(\d\sigma(t)\right)^2\le 0,
$$
for all $u\in U$, almost every $t\in [0,T]$ and $\mathbb{P}-$almost surely.

\ms\no
This finishes the proof of Theorem \ref{main result}.   $\qed$

\section{Illustrative Example:  Linear-quadratic risk-sensitive model }  \label{sec:example}
The optimal control of a linear stochastic system driven by a Brownian
motion and with a quadratic cost in the state and the control is probably
the most well known solvable stochastic control problem in continuous
time. To illustrate our approach, we consider the one-dimensional case with linear state dynamics and exponential quadratic cost functional.

\ms  It is well-known that in absence of mean-field coupling, the optimal control is a linear feedback control whose feedback gain is obtained from the solution of a risk-sensitive Riccati equation which has an additional term when compared to the (classical) Riccati equation for the quadratic cost problem. In the examples below we will show that this feature is still valid in the LQ risk-sensitive problem (with and without the mean-field coupling).

\subsection{LQ risk-sensitive control without the mean-field coupling}

We consider the linear-quadratic risk-sensitive control problem:
\begin{eqnarray}
\label{LQ} 
\left\{
\begin{array}{lll} \inf_{u(\cdot)\in \cU}  E e^{\theta \left[ \frac{1}{2} \int_0^T u^2(t)dt +\frac{1}{2}x^2(T)\right]},
\\
\displaystyle{\mbox{ subject to }\ }\\
dx(t)=\left(ax(t)+bu(t)\right)dt+\sigma dB_t,\\
x(0)=x_{0},\\
 \end{array}
\right.
\end{eqnarray}
where, $a, b$ and $\sigma$ are real constants.

\ms\no An admissible pair $(\bar x(\cdot), \bar u(\cdot))$ that satisfies the optimality  necessary conditions of Theorem \ref{main result} can be obtained by solving the following system of forward-backward SDEs: 

\be\label{LQ-AD}
\left\{ \begin{array}{lll}
d\bar x(t)=\left(a\bar x(t)+b\bar u(t)\right)dt+\sigma dB_t,\\
dp(t)=-\left\{ ap(t)+\theta\ell(t)q(t)\right\}dt+q(t)dB_t,\\
dv^{\theta}(t)=\theta\ell(t)v^{\theta}(t)dB_t, \\
v^{\theta}(T)=\phi^{\theta}(T),\\
x(0)=x_0,\quad p(T)=-\bar x(T),
\end{array}
\right.
\ee
where, $\phi^{\theta}(T):=e^{\theta \left[ \frac{1}{2} \int_0^T \bar u^2(t)dt +\frac{1}{2}\bar x^2(T)\right]}$.

\ms \no This system  involves only the first adjoint equation because the diffusion coefficient in the state dynamics  is independent of the control (constant!).

\medskip The associated risk-sensitive Hamiltonian is 
$$
H^{\theta}(t,x,u, p, q, \ell):=(ax+bu)p- \frac{1}{2}  u^2+\sigma(q+\theta\ell p).
$$
We have 
$$
H^{\theta}_x=ap, \qquad  H^{\theta}_u=bp-u.
$$
Maximizing the Hamiltonian yields 
\be\label{LQ-u}
\bar u(t)=bp(t).
\ee
The associated state dynamics $\bar x$ solves then the  SDE
\be\label{LQ-state}
d\bar x(t)=\left(a\bar x(t)+b^2p(t)\right)dt+\sigma dB_t.
\ee

\ms\no  We try a solution of the form 
\be\label{LQ-p}
p(t):=-\b(t)\bar x(t), 
\ee
where, $\beta(t)$ is a deterministic function such that $\beta(T)=1$. In view of (\ref{LQ-p}), the state dynamics $\bar x$ solves the linear SDE
\be\label{linear}
d\bar x(t)=\left(a-b^2\b(t)\right)\bar x(t)dt+\sigma dB_t.
\ee
Furthermore,  we have
\be\label{LQ-p-x}
dp(t)=\left(-\dot\b(t)-a\b(t)+b^2\b^2(t)\right)\bar x(t)dt-\sigma\beta(t)dB_t.
\ee
From (\ref{LQ-AD}) we also get
\be
dp(t)=-\left\{-a\b(t)\bar x(t)+\theta\ell(t)q(t)\right\}dt+q(t)dB_t.
\ee
Identifying the coefficients of these two equations we obtain
\be\label{LQ-q}
q(t)=-\sigma\b(t),
\ee
and
$$
\left(-\dot\b(t)-2a\b(t)+b^2\b^2(t)\right)\bar x(t)=\sigma\b(t)\theta\ell(t).
$$
This equation is feasible only if we choose 
\be\label{LQ-ell}
\ell(t)=\gamma(t)\bar x(t).
\ee
for some deterministic function $\gamma(t)$. 

\medskip\noindent  Given the deterministic functions $\b$ and $\gamma$, in view of (\ref{LQ-AD}), (\ref{linear}) and (\ref{LQ-ell}), the generic martingale $v^{\theta}$ satisfies the linear SDE
\be\label{LQ-v}
dv^{\theta}(t)=\theta\gamma(t)\bar x(t)v^{\theta}(t)dB_t,\quad v^{\theta}(T)=\phi^{\theta}(T).
\ee 
At this stage, $\gamma$ can be seen as a free parameter whose choice gives different features of the behavior of the optimal pairs $(\bar x(\cdot),\bar u(\cdot))$. Let us examine two typical cases (among many others).

\medskip {\bf Case 1}.  $\gamma(t):=\sigma\b(t)$. 

\noindent  This choice yields the form  $\ell (t)=-\sigma p(t)$ suggested by Lim and Zhou \cite{lim}, using the relationship between the SMP and the dynamic programming principle, which in turn gives the  risk-sensitive  Riccati equation for $\b$:
$$
\dot\b(t)+2a\b(t)+(\theta\sigma^2-b^2)\b^2(t)=0,\qquad \b(T)=1.
$$ 
Its  explicit solution is given by 
$$
\b(t)=\left[\frac{b^2-\theta\sigma^2 }{2a}+\left(1-\frac{b^2-\theta\sigma^2 }{2a}\right)e^{-2a(T-t)}\right]^{-1}.
$$

\medskip {\bf Case 2}.  $\gamma(t):=1$. 

\noindent  This choice yields the form  $\ell(t)=\bar x(t)$, which is not related to the choice made in Lim and Zhou \cite{lim}. We obtain  yet another Ricatti equation for $\b$:
$$
\dot\b(t)+(2a+\theta\sigma)\b(t)-b^2\b^2(t)=0,\qquad \b(T)=1.
$$
Its explicit solution is given by
$$
\b(t)=e^{(2a+\theta\sigma)(T-t)}\left[1-\frac{b^2e^{\theta\sigma T}}{2a+\theta\sigma}\left(e^{(2a+\theta\sigma)(T-t)}-1\right)\right]^{-1}.
$$

\no Note that, depending on the parameters, there may be explosion of the function $\b$ in finite time, in both solutions.

\subsection{LQ risk-sensitive control with a mean-field coupling}
We keep the same functions $f,\sigma,b$ as in (\ref{LQ})  but we modify the terminal cost $h$ to be 
$$
h(x(T), E[x(T)])= \frac{1}{2}x^2(T)+\mu E[x(T)],
$$
for some given constant $\mu$,  where the only mean field coupling is $E[x(T)]$.  Therefore, the first-order adjoint equation remains the same as (\ref{LQ-AD}), but  the terminal condition  becomes
\be\label{T-p}
p(T)= -\bar x(T)-\frac{\mu}{\phi^{\theta}_T}E[\phi^{\theta}_T],
\ee
where,  $\phi^{\theta}_T:=e^{\theta \left[ \frac{1}{2} \int_0^T \bar u^2(t)dt +\frac{1}{2}\bar x^2(T)+\mu E[\bar x(T)]\right]}$.

\medskip\noindent In view of (\ref{mart-rep}), we have
\be\label{exp-L}
L^{\theta }_t:=\frac{v^{\theta}(t)}{E[\phi^{\theta}_T]}=\frac{v^{\theta}(t)}{v^{\theta}(0)}=\exp{\left(\theta\int_0^t\ell(s)dB_s-\frac{\theta^2}{2}\int_0^t |\ell(s)|^2ds\right)}.
\ee
 which satisfies the linear SDE
\be\label{L-SDE}
dL^{\theta }_t=\theta \ell(t) L^{\theta }_tdB_t, \qquad L^{\theta }_0=1.
\ee
Hence, the end-value (\ref{T-p}), becomes
\be\label{T-p-1}
L^{\theta }_Tp(T)=-L^{\theta }_T\bar x(T)-\mu.
\ee

\medskip\noindent The associated risk-sensitive Hamiltonian is 
$$
H^{\theta}(t,x,u, p, q, \ell):=(ax+bu)p- \frac{1}{2}  u^2+\sigma(q+\theta\ell p).
$$
We have 
$$
H^{\theta}_x=ap, \qquad  H^{\theta}_u=bp-u.
$$
Maximizing the Hamiltonian yields 
\be\label{MF-LQ-u}
\bar u(t)=bp(t).
\ee

\ms\no In view of the form (\ref{T-p-1}), we try a solution $p(t)$  such that
\be\label{MF-LQ-p}
p(t) L^{\theta}_t= - \b(t) \bar x(t)L^{\theta}_t- \mu \a(t),
\ee
where, $\a(t)$ and $\b(t)$ are deterministic function such that  $\a(T)= 1$ and $\b(T)=1$.  Proceeding as above, in view of (\ref{MF-LQ-u}), we apply It\^o's formula to the process $v^{\theta}(t)p(t)$ using  (\ref{LQ-AD}) and (\ref{L-SDE}), and  then ({\ref{MF-LQ-p}}) and (\ref{LQ-state}),  and identify the coefficient we obtain
$$\left\{\begin{array}{lll}
q(t)=\frac{\theta \ell(t)}{ L^{\theta}_t}\mu\a(t)-\sigma\b(t),\\
-\mu\dot\a(t)+\mu(b^2\b(t)-a)\a(t)-\left[\left(\dot\b(t)+2a\b(t)-b^2\b^2(t)\right)\bar x(t)+\sigma\theta\ell(t)\b(t)\right]L^{\theta}_t=0.
\end{array}
\right.
$$
Again, as in the previous example,  this system of equations is feasible only if we assume that
\be\label{MF-LQ-ell}
\ell(t)=\gamma(t)\bar x(t),
\ee 
for some deterministic function $\gamma(t)$. This yields
\be\label{MF-LQ-diff}\left\{\begin{array}{lll}
\dot\b(t)+2a\b(t)-b^2\b^2(t)+\theta\sigma\gamma(t)\b(t)=0,\\
\dot\a(t)+(a-b^2\b(t))\a(t)=0,\\
q(t)=\mu \theta\gamma(t)\bar x(t) (L^{\theta}(t))^{-1}\a(t)-\sigma\b(t),\\
\b(T)=1, \quad \a(T)=1.
\end{array}
\right.
\ee
Finally, choosing either $\gamma(t)=\sigma\b(t)$ or $\gamma(t)=1$, as in the previous example, we get closed form solutions to the system (\ref{MF-LQ-diff}).

\medskip\no We note that when $\mu$ goes to zero, we obtain the mean-field free solution given in the previous subsection. Moreover, the choice of the process $\ell$ need not be related to any relationship between the stochastic maximum principle and the Dynamic Programming Principle.

\section*{Numerical investigation}
\begin{figure}[htb]
	\centering
		\includegraphics[scale=1.2]{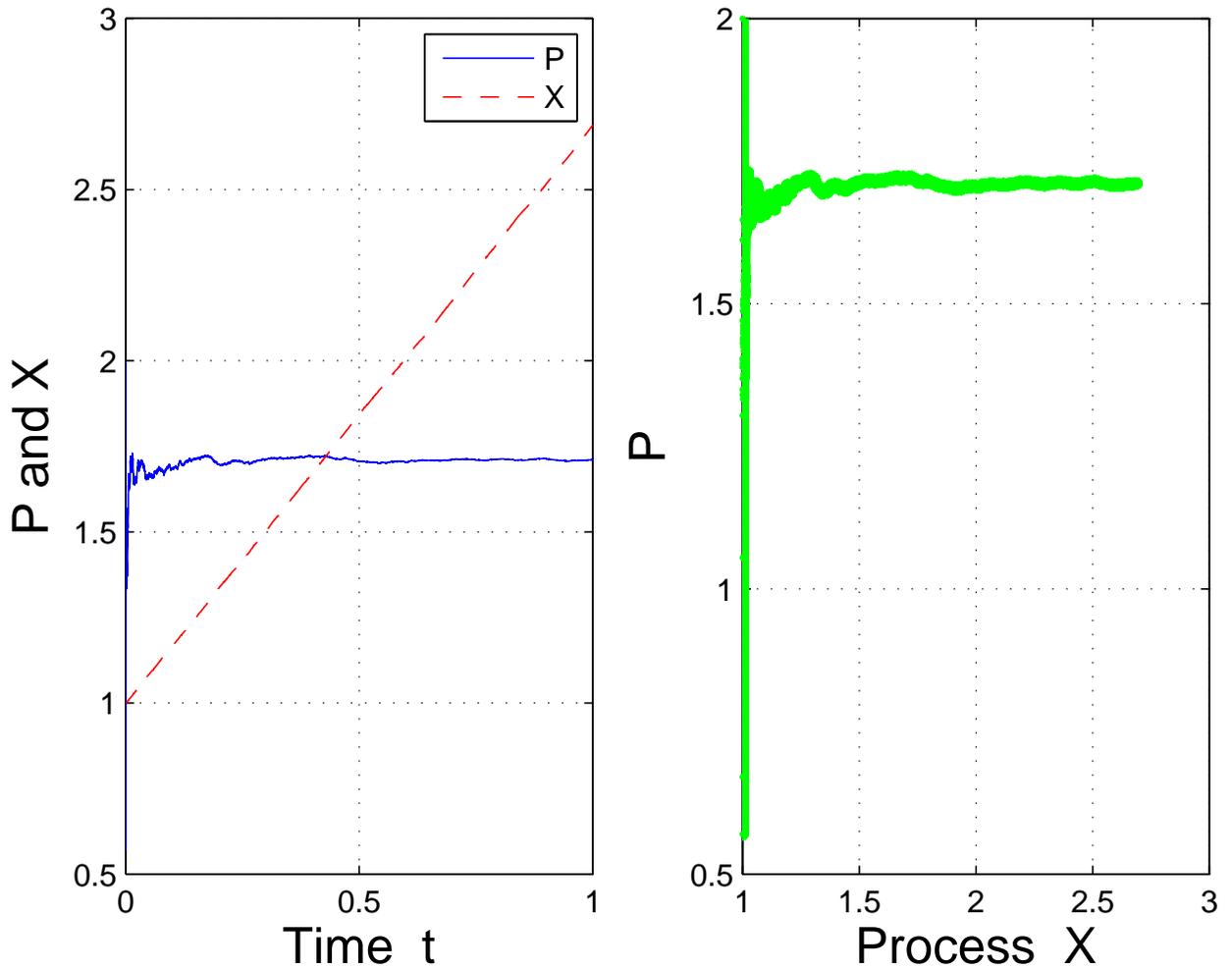}
	\caption{A local solution to the risk-sensitive linear quadratic system with small parameter $\theta.$}
	\label{figure1}
\end{figure}

\begin{figure}[htb]
	\centering
		\includegraphics[scale=1.2]{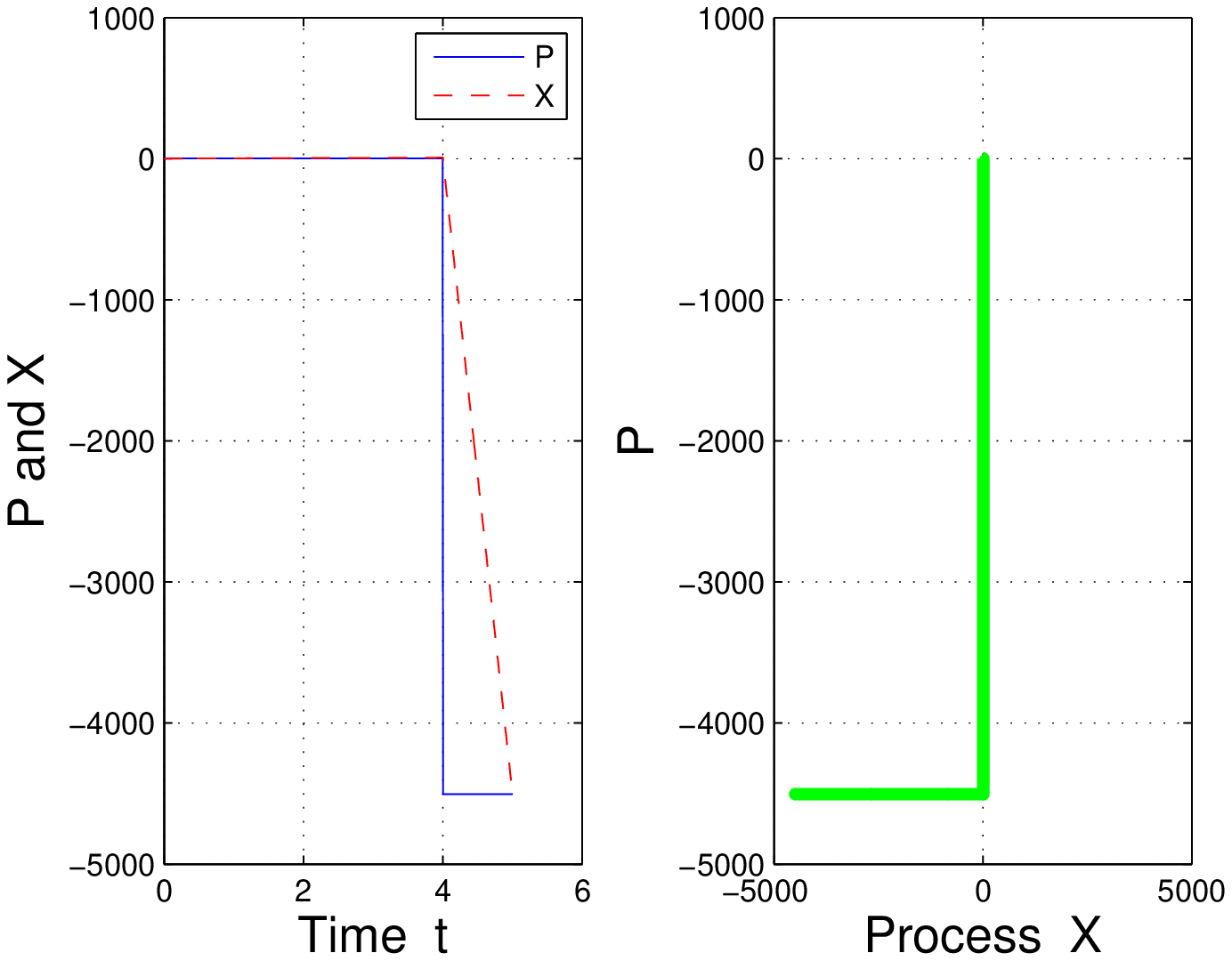}
	\caption{Explosion of the local solution to the risk-sensitive linear quadratic system when $T$ is larger.}
	\label{figure2}
\end{figure}

In this subsection we provide numerical solution of the above risk-sensitive linear quadratic system. The state parameters are $a=0, b=1.$ The initial state is $1.$ The noise parameter is   $\sigma=10^{-2}$ and the risk-sensitive index is set to $\theta=10^{-5}.$ The step size of the discretization is set to $10^{-6}.$
We observe that a local solution exist for small window $[0, 1]$ as illustrated in Figure \ref{figure1}. As $T$ becomes larger (for example for $T=5$ in Figure \ref{figure2} ) there is an expolosion of the solution.

\section{Conclusion} \label{sec:conclusion}
In this paper we established a Peng's type stochastic maximum principle for risk-sensitive stochastic control of mean-field type extending a previous result by Lim and Zhou \cite{lim}.

\bibliographystyle{plain}


\end{document}